\documentclass[12pt,technote, onecolumn]{IEEEtran}

\usepackage{amssymb}
\usepackage{theorem}
\usepackage{amsmath} 
\usepackage{listings}
\usepackage{multicol}
\usepackage{array}
\usepackage{calc}
\usepackage{ifthen}
\usepackage{url}
\usepackage{amsmath}
\pagebreak[5]
\title{Isometries and Construction of Permutation Arrays} \author{Mathieu Bogaerts}
\usepackage{algorithmic}
\usepackage{algorithm}
\newtheorem{theorem}{Theorem}
\begin{document}
\maketitle
\begin{abstract}
An $(n,d)$-permutation code  is a subset $C$ of $Sym(n)$ such that the Hamming distance $d_H$ between any two distinct elements of $C$ is at least equal to $d$. In this paper, we use the characterisation of the isometry group of the metric space $(Sym(n),d_H)$ in order to develop generating algorithms with rejection of isomorphic objects. To classify the $(n,d)$-permutation codes up to isometry, we construct invariants and study their efficiency. We give the numbers of non-isometric $(4,3)$- and $(5,4)$- permutation codes. Maximal and balanced $(n,d)$-permutation codes are enumerated in a constructive way. 
\end{abstract}
\section{Permutation arrays and permutation codes}
 
An $(n,d)$-permutation code of distance $d$, size $s$ and degree $n$ is a non-empty subset $C$ of the symmetric group $Sym(n)$ such that the Hamming distance between any two distinct elements of $C$ is at least equal to $d$. The Hamming distance between two permutations $\phi,\psi\in Sym(n)$ is defined as
 $d_H(\phi, \psi)=|\{ i \in \{1,\dots,n\}: \phi(i) \ne\psi(i)  \}| $.
The $s\times n$ array $A$ associated to a $(n,d)$-permutation code $C=\{\phi_1,\dots,\phi_s\}$ of size $s$ by $A_{ij}=\phi_i(j)$ has the following properties: every symbol $1$ to $n$ occurs exactly in one cell of any row and any two rows disagree in at least $d$ columns.
Such an array is called a permutation array (PA) of distance $d$, size $s$ and degree $n$.
Permutation codes have first been proposed by Ian Blake in 1974 as  error-correcting codes for powerline communication \cite{D29}. This application motivates the study of permutation codes having the largest possible size. Upper bounds for the maximal size $\mu(n,d)$ of a permutation code with fixed parameters $n$ and $d$ have been studied by many authors, see e.g. Deza \cite{C26}, Cameron \cite{C3}, and more intensively since Chu, Colbourn and Dukes \cite{C6},Tarnanen \cite{C14}, and Han Vinck \cite{B19, B22}.\\

A $(n, d)$-permutation code $C$ of size $s$ is maximal if $C$ is not contained in a $(n, d)$-permutation code of larger size $s'>s$. Note that a $(n, d)$-permutation code reaching the maximal size $\mu(n,d)$ is necessarily maximal while the converse is not true. In this article, we will focus on the generation up to isometry of all maximal permutation codes with fixed parameters $n$ and $d$. \\

The paper is organised as follows: section \ref{secIso} describes the isometry group of the metric space $(Sym(n), d_H)$, section \ref{secGenAllPA} contains improvements of methods to obtain all PA's of a given size. Section \ref{GenAllMaxPA} describes two algorithms for generating all maximal permutation codes up to isometry. 
\section{Isometries}\label{secIso}

A distance $D$ on $Sym(n)$ is called left-invariant (resp. right-invariant) if $D(\phi, \psi)=D(\alpha \phi,\alpha \psi)$ (resp. $D(\phi, \psi)=D(\phi\alpha,\psi\alpha)$ ) for all $\alpha,\phi,\psi \in Sym(n)$. A distance that is both left- and right-invariant is said to be bi-invariant. For any bi-invariant distance, the left multiplications $l_\alpha:\phi\mapsto \alpha\phi$  and the right multiplications $r_\alpha: \phi\mapsto \phi\alpha^{-1}$ are isometries.  As noticed by Deza \cite{D18}, any bi-invariant distance is  invertible: $D (\phi,\psi)= D(\phi^{-1} ,\psi^{-1})$, or equivalently, the inversion $\i$, mapping each permutation onto its inverse, is an isometry. Let denote by $\mathcal{R}$ (resp. $\mathcal{L}$) the group of all right (resp. left-) multiplications and by $\mathcal{I}$ the group generated by the inversion $\i$. We will say that the distance $D$  distinguishes the transpositions if there exists a constant $c$ such that  $D(\phi,\psi)= c \iff \phi\psi^{-1}$ is a transposition. \\

In 1960, Farahat characterized  the isometry group $Iso(n)$ of the metric space $(Sym(n),d_H)$  \cite{D17}. Since the Hamming distance is  bi-invariant and distinguishes the transpositions, the following result generalizes the characterisation given by Farahat:
\begin{theorem}\label{TheoIso}
Let $D$ be a bi-invariant  distance distinguishing the transpositions on $Sym(n)$ ($n\ge 3$), then the group $Iso_D$ of isometries of $(Sym(n), D)$ is $(\mathcal{L} \times \mathcal{R})\rtimes \mathcal{I}$, isomorphic to $Sym(n) \wr 2$.
\end{theorem}

Proof: For any bi-invariant  distance $D$ distinguishing the transpositions, the isometry group of $(Sym(n),D)$  contains the groups $\mathcal{L}$, $\mathcal{R}$ and $\mathcal{I}$, generating the group $(\mathcal{L} \times \mathcal{R})\rtimes \mathcal{I}$. The pointwise stabiliser of the set containing the identity, the transpositions $(1,k)$ with $2\le k\le n$ and $(1,2,3)$ in the isometry group $Iso(n)$ is the trivial group, proving that there is no other isometries than $(\mathcal{L} \times \mathcal{R})\rtimes \mathcal{I}$. See \cite{MB2009} for developements.

Every isometry $t\in Iso(n)$ can be uniquely written as $l_\alpha r_\beta \i^k$ with $k=0$ or $1$, $\alpha,\beta\in Sym(n)$. The action of a left multiplication $l_\alpha$ on a given code corresponds to the permutation under $\alpha$  of the symbols appearing in the PA associated to the code, and the action of a right-multiplication $r_\beta$ is equivalent to the permutation under $\beta$ of the columns of the PA. In other words, classifying  permutation codes up to isometry is equivalent to classifying PA's up to permutation of their rows, their columns, their symbols and up to the inversion. \\ 
\section{Generate all $(n, d)$-permutation codes  of a given size}\label{secGenAllPA}

An efficient way to generate all permutation codes of a given size is described in \cite{C6} as a clique search. A graph is built whose vertices are the elements of $Sym(n)$ and whose edges are the pairs of permutations at distance at least equal to $d$. A (maximal) $(n,d)-$permutation code is equivalent to a (maximal) clique in this graph.  Since the isometry group $Iso(n)$ is transitive on $Sym(n)$, any $(n, d)$-permutation code is isometric to a permutation code containing the identity permutation $Id$.  Assuming that the permutation code contains $Id$, the clique search is restricted  to the neighbourhood of the vertex $Id$, that is, to the set of permutations that are at distance at least equal to $d$ from the permutation $Id$. This improvement reduces the number of vertices to consider from $n!$ to $\sum_{k=d}^n {n \choose k}D_k$, where $D_k$ is the number of fixed-point-free permutations of $k$ points .\\

Using this method, we confirm the results of T. Kl{\o}ve \cite{B5}:  there are, up to isometry, exactly 7 $(6,5)-$ permutation codes of size 18. Note that the equivalence between PA's used in \cite{B5} doesn't take in consideration the inversion $\i$.\\

The subgroups of $Sym(n)$  can be used to construct large PA's. Starting form a subgroup $G$, Chu, Colbourn and Dukes used the orbits of $L_G=<l_\beta : \beta \in G>$ to simplify the construction of PA's using a clique search algorithm. More recently, Janiszczak and Staszewski \cite{IJRS} used the orbits a subgroup of $Iso(10)$ (more precisely, a subgroup $U_G$ of the conjugacies $U_G=<l_\beta r_{\beta}: \beta \in G>$) in order to prove that $\mu(10,9)\ge 49$. Using other subgroups of $Iso(n)$, we found other codes as e.g. an $(7,4)-$permutation code of size $336$ preserved by a subgroup of $Iso(7)$ of order $7056$ \cite{HomepageMB}.

\section{Generating all maximal permutation codes up to isometry}\label{GenAllMaxPA}

In the past few years, techniques have been developped in order to construct PA's of large size. We adopt another approach, and for small values of $n$ and $d$, succeed to classify all maximal permutation codes up to isometry.
Note that for $d=n$, a maximal permutation code has size $s=n$ and the PA associated to such a code is a latin square. The enumeration of latin squares up to permutation of the rows, columns and symbols has been widely  studied and has been achieved for $n\le 11$, see \cite{B16, B14}. As far as we know, the only classification of maximal $(n,d)$-permutation codes with $d<n$ has been made by Torleiv Kl{\o}ve for $(6,5)-$permutation codes of size $18$, as mentionned in section \ref{secGenAllPA}. \\

In order to generate all maximal permutation codes with fixed parameters $n$ and $d$,  we used two isomorphism free generation algorithms, as introduced by McKay \cite{D13}. The first one is based on the backtrack principle, at every step the generated code is compared to the list of already encountered codes. We introduce here the notation $V_d(C)$ for the set of permutations that are at distance at least equal to $d$ of all elements in $C$. \\  
\begin{algorithm}
\caption{Generate a list $L$ of all non-isometric codes}
\begin{algorithmic}\label{GENBYLIST} 
\STATE { Procedure} GENBYLIST ($C, L$)
\IF{ $\exists C'  \in L$ such that $C\sim C'$ } 
\STATE break
\ENDIF 
\STATE $L \leftarrow L \cup C$
\IF{ $C$ is maximal}
\STATE return $C$
\ENDIF
\FOR{$\phi \in V_d(C)$}
\STATE GENBYLIST( $ C\cup \{\phi\},L$)
\ENDFOR
\end{algorithmic}
\end{algorithm}

Initiated with an empty list and  a code $C$ consisting of a single permutation, algorithm \ref{GENBYLIST} produces exactly one representative of each permutation code up to isometry. It requires to compare permutation codes,  motivating the developpement of invariants that are efficient to test whether two permutation codes are isometric or not.\\

We define the pair of quotient sets of a $(n,d)$-permutation code $C$ as the ordered set  $(\Delta(C),\Sigma(C))$ where $$\Delta(C):=\{\phi\psi^{-1}: \phi,\psi \in C\}$$
$$\Sigma(C):=\{\phi^{-1}\psi: \phi,\psi \in C\}$$
Two pairs of quotients sets $(\Delta, \Sigma ),(\Delta', \Sigma' )$ are equivalent if there exists $\alpha, \beta\in Sym(n)$ such that $\alpha \Delta \alpha^{-1}=\Delta'$ and $\beta \Sigma \beta^{-1}=\Sigma'$ or $\alpha \Delta \alpha^{-1}=\Sigma'$ and $\beta \Sigma \beta^{-1}=\Delta'$.

A famous code invariant is the distance enumerator polynomial of a code $C$, that enumerates the number of pairs of permutations wich are at distance $k$ for $0\le k\le n$. This invariant is weaker than the cycle index polynomial $Q_C(x)$ of the permutation code $C$, defined using the $p(n)$ conjugacy classes $C_1,\dots, C_{p(n)}$ as  $Q_C(x)=\frac{\sum_{j=1}^{p(n)} b_j x^j}{|C|}$ where $b_j$ is the number of pairs of permutations $(\phi,\psi)$ of $C$ such that $\phi\psi^ {-1}\in C_j$. If the code $C$ is a group, the cycle index polynomial of $C$ coincides with the usual definition of the cycle index polynomial of a group.  \\

The set of occurences $Occ(C)$ associated to a $(n,d)$-permutation code $C$ is defined as $Occ(C):=\{o_{ij}:1\le i, j\le n\} $ where $o_{ij}$ is the number of permutations $\phi\in C$ such that $\phi(i)=j$. If $Occ(C)=\{r\}$, then the permutation code is said to be $r$-balanced (see \cite{C24}). \\

Given an $(n,d)$-permutation code $C$ of size $s$, we define the graph $(V(C),E(C))$ as follows. The $s+n^2+2n$ vertices are partionned into four sets $V_1:=\{\phi_i: 1\le i\le s\}$, $V_2:=\{c_j: 1\le j\le n\}$, $V_3:=\{ s_j: 1\le j\le n\}$ and $V_4:=\{c_{ij}: 1\le i, j\le n  \}$, $\phi_i$ is  the $i$-th permutation of $C$, $c_j$ is the $j$-th column, $s_j$ is the $j$-th symbol and $c_{ij}$ is the $(i,j)$-cell of the permutation array associated to $C$. The edges of the graph are defined by the following rule:  $\phi_k(i)=j$ if and only if the vertex $c_{ij}\in V_4$ is adjacent to the vertices $\phi_k\in V_1$, $c_i\in V_2$ and $s_j\in V_3$.

\begin{theorem}
The pair of quotient sets, the cycle index polynomial and the set of occurences are invariants, but are not complete invariants. The graph $(V(C), E(C))$ of the code $C$ is a complete invariant.
\end{theorem}
The proof of this statement is routine check since the isometry group is described by Theorem \ref{TheoIso}.

The graph $( V(C),E(C))$ can also be used to obtain a canonical form of the permutation code by a canonical labelling of the vertices of the graph (using e.g. routines of Nauty \cite{Nauty}). The canonical form automatically defines a canonical augmentation $p(C)$ for any code $C$, and insures that the algorithm \ref{CANAUG} generates all non-isomorphic maximal codes.\\
\begin{algorithm}
\caption{Returns all non-isometric maximal codes}
\begin{algorithmic}\label{CANAUG} 

\STATE { Procedure} CANONICAL-AUGMENTATION ($C$)
\IF{ $C$ is maximal}
\STATE return $C$
\ENDIF
\STATE $\mathcal{C}  \leftarrow \{ C\cup\{\phi\} :\phi \in V_d(C) \}$
\STATE $\mathcal{D} \leftarrow \{\mathcal{C} \cap \{ t(D): t\in Stab(C) \} : D\in \mathcal{C} \}$
\FOR{ $Z \in \mathcal{D}$}
\STATE choose $K\in Z$
\IF{ $\exists t\in Stab(K)$ s. t. $t(C)= p(C)$ }
\STATE CANONICAL-AUGMENTATION($K$)
\ENDIF
\ENDFOR
\end{algorithmic}
\end{algorithm}\\
Results of algorithm \ref{GENBYLIST} are confirmed by the results of algorithm \ref{CANAUG}: there are, up to isometry, $61$ $(4,3)-$permutation codes ($4$ are maximal) and $9445$ $(5,4)-$permutation codes ($139$ are maximal).
Table \ref{NonIsoPA} details the numbers of non-isometric maximal PA's of distance 4 on $Sym(5)$, of size 7 to 20.
\begin{table}[!h]
\centering
\caption{Non-isometric maximal $(5,4)-$PA's }\label{NonIsoPA}
\begin{IEEEeqnarraybox}[\IEEEeqnarraystrutmode\IEEEeqnarraystrutsizeadd{2pt}{1pt}]{v/c/v/c/v}
\IEEEeqnarrayrulerow\\
&\mbox{Size}&&\mbox{Number of non-isometric $(5,4)-$PA's}&\\
\IEEEeqnarraydblrulerow\\
& 7 && 1&\\
&8&& 25&\\
&9&& 36&\\
&10&&46&\\
&11&&18&\\
&12&& 10&\\
&13&&1&\\
&15&&1&\\
&20&&1&\\
\IEEEeqnarrayrulerow
\end{IEEEeqnarraybox}
\end{table}

An $r$-balanced PA with $n$ symbols has size exactly $nr$. We constructed all non-isometric $r$-balanced codes $C(n, d)$ for small values of $n$ and $d$. Table \ref{NonIsoPAbal} shows the exact number $M_r(n,d)$ of non-isometric $r-$balanced $(n, d)-$permutation codes for small values of $n, d$ and $r$. All our permutation codes and the code of our Magma \cite{MAGMA} programs are available \cite{HomepageMB}.
\begin{table}[h!t]
\centering
\caption{Number of non-isometric $r$-balanced PA's}\label{NonIsoPAbal}
\begin{IEEEeqnarraybox}[\IEEEeqnarraystrutmode\IEEEeqnarraystrutsizeadd{2pt}{1pt}]{v/c/v/c/v/c/v/c/v}
\IEEEeqnarrayrulerow\\
& r && n && d && M_r(n,d) &\\
\IEEEeqnarraydblrulerow\\
&$2$ && $5$ && $4$ && $6$&\\
&$2$ && $5$ && $3$ && $218$&\\
&$2$ && $6$ && $5$ && $2799$&\\
&$3$ && $5$ && $4$ && $1$&\\
&$4$ && $5$ && $4$ && $1$&\\
\IEEEeqnarrayrulerow
\end{IEEEeqnarraybox}
\end{table}\\

\providecommand{\bysame}{\leavevmode\hbox to3em{\hrulefill}\thinspace}
\providecommand{\MR}{\relax\ifhmode\unskip\space\fi MR }
\providecommand{\MRhref}[2]{%
  \href{http://www.ams.org/mathscinet-getitem?mr=#1}{#2}
}
\providecommand{\href}[2]{#2}


\begin{thebibliography}{10}

\bibitem{B19}
V.~B. Balakirsky and A.~J.~Han Vinck, \emph{On the performance of permutation
  codes for multi-user communication}, Probl. Inf. Transm. \textbf{39} (2003),
  no.~3, 239--254.

\bibitem{D29}
I.~Blake, \emph{Permutation codes for discrete channels}, IEEE Tansactions on
  Information Theory (1974), 138--140.

\bibitem{HomepageMB}
M.~Bogaerts, \emph{Homepage}, \url{http://homepages.ulb.ac.be/~mbogaert/}.

\bibitem{MB2009}
\bysame, \emph{Codes et tableaux de permutations: construction ,
  \'enum\'eration et automorphismes}, Ph.D. thesis, Universit\'e Libre de
  Bruxelles, June 2009,
  \url{http://theses.ulb.ac.be/ETD-db/collection/available/ULBetd-06252009-090%
241/}.
\bibitem{MAGMA}
W. ~Bosma, J. ~Cannon, and C. ~Playoust, \emph{The Magma algebra system. I. The user language.}, J. Symbolic Comput. \textbf{24},3-4 (1997), 235--265
\bibitem{C3}
P.~J. Cameron, \emph{Metric and geometric properties of sets of permutations},
  Algebraic, Extremal and Metric Combinatorics 1986 (Deza, Frankl, and
  Rosenberg, eds.), Cambridge University Press, 1988, pp.~39--53.

\bibitem{C6}
W.~Chu, C.~J. Colbourn, and P.~Dukes, \emph{Construction for permutation codes
  in powerline communications}, Des. Codes Cryptography \textbf{32} (2004),
  51--64.

\bibitem{C26}
M.~Deza and P.~Frankl, \emph{On the maximum number of permutations with given
  maximal or minimal distance}, Journal of Combinatorial Theory (A) \textbf{22}
  (1977), 352--360.

\bibitem{D18}
M.~Deza and T.~Huang, \emph{Metrics on permutations, a survey}, J.
  Combinatorics, Information and System Sciences \textbf{23} (1998), 173--185.

\bibitem{C24}
C.~Ding, F.-W. Fu, T.~Kl{\o}ve, and V.~K.-W. Wei, \emph{Construction of
  permutation arrays}, IEEE Transactions on Information Theory \textbf{48}
  (2002), no.~4, 977--980.

\bibitem{D17}
H.~Farahat, \emph{The symmetric group as a metric space}, Journal of London
  Mathematical Society \textbf{35} (1960), 215--220.

\bibitem{IJRS}
I.~Janiszczak and R.~Staszewski, \emph{An improved bound for permutation arrays
  of length 10}, \url{http://www.iem.uni-due.de/preprints/IJRS.pdf}.

\bibitem{B5}
T.~Kl{\o}ve, \emph{Classification of permutation codes of length 6 and minimum
  distance 5}, Proceedings of The 2000 International Symposium on Information
  Theory and Its Application (ISITA '2000), 2000, pp.~465--468.

\bibitem{D13}
B.~D. McKay, \emph{Isomorphism-free exhaustive generation}, Journal of
  algorithms \textbf{26} (1998), 306--324.

\bibitem{Nauty}
\bysame, \emph{Nauty (version 2.2)}, 2008, Software
  \url{http://cs.anu.edu.au/~bdm/nauty}.

\bibitem{B16}
B.~D. McKay, A.~Meynert, and W.~Myrvold, \emph{Small latin squares, quasigroups
  and loops}, Journal of Combinatorial Designs \textbf{15} (2007), no.~2,
  98--119.

\bibitem{B14}
B.~D. McKay and I.~M. Wanless, \emph{On the number of latin squares}, Annals of
  Combinatorics \textbf{9} (2005), no.~3, 335--344.

\bibitem{C14}
H.~Tarnanen, \emph{Upper bounds on permutation codes via linear programming},
  European J. Combinatorics \textbf{20} (1999), 101--114.

\bibitem{B22}
A.J.~Han Vinck, \emph{Coded modulation for power-line communications}, AEÜ Int.
  J. Electron. and Commun. \textbf{54} (2000), no.~1, 45--49.

\end{thebibliography}
\end{document}